\documentclass[11pt]{elsarticle}
\usepackage{amsfonts,amssymb,amsmath}
\usepackage{graphicx}
\usepackage{psfrag}
\usepackage{calc}
\usepackage{epic}

 \usepackage{caption2}

\usepackage{comment}

\newtheorem{theorem}{Theorem}
\newtheorem{lemma}[theorem]{Lemma}

\newtheorem{proposition}[theorem]{Proposition}

\newtheorem{assumption}{Assumption}

\def\R{ {\mathbb{R} }}

\def\N{ {\mathbb{N} }}

\def\E{\mathbb{E}}
\def\1{\mathbb{1}}

\def\clL{ {\mathcal{L} } }

\def\psid{ {\psi^\circ}}

\def\thetad{ {\theta^{\downarrow}} }
\def\thetau{ {\theta^{\uparrow}} }

\def\phisd{ {\phi^{\downarrow} }}
\def\phisu{ {\phi^{\uparrow} }}

\def\phidsd{ {\dot\phi^{\downarrow} }}
\def\phidsu{ {\dot\phi^{\uparrow} }}

\def\gammasd{ {\gamma^{\downarrow} }}
\def\gammasu{ {\gamma^{\uparrow} }}

\def\psisd{ {\psi^{\downarrow} }}
\def\psisu{ {\psi^{\uparrow} }}

\def\wbar{ {\bar{w}} }

\def\pf#1{\paragraph{\rm\scshape#1}}
\def\sq{\hbox{\rlap{$\sqcap$}$\sqcup$}}
\def\qed{\ifmmode\sq\else{\unskip\nobreak\hfil
\penalty50\hskip1em\null\nobreak\hfil\sq
\parfillskip=0pt\finalhyphendemerits=0\endgraf}\fi\medskip}

\def\epsy{\varepsilon}
\def\varble{\,\cdot\,}
\input{epsf}

\def\EboxRatePath#1#2{%
\begin{center}
\parbox[t]{.35\hsize}{\epsfxsize=\hsize \epsfbox{#1}}
\qquad\qquad
\parbox[t]{.38\hsize}{\epsfxsize=\hsize \epsfbox{#2}}
\end{center}}

\def\Ebox#1#2{%
\begin{center}
\parbox{#1\hsize}{\epsfxsize=\hsize \epsfbox{#2}}
\end{center}}

 \newlength{\noteWidth}
\long\def\notes#1{\ifinner
             {\tiny #1}
             \else
              \marginpar{\parbox[t]{\noteWidth}{\raggedright\tiny #1}}
               \fi}
\def\notes#1{\typeout{#1 !!!}}

\def\Prob{{\sf P}}

\def\Expect{{\sf E}}

\def\eqdef{\mathbin{:=}}

\def\Re{\mathbb{R}}

\def\Lemma#1{Lemma~\ref{thm:#1}}
\def\Theorem#1{Theorem~\ref{thm:#1}}
\def\Proposition#1{Proposition~\ref{thm:#1}}

\def\Figure#1{Figure~\ref{fig:#1}}

 \def\FRAC#1#2#3{\genfrac{}{}{}{#1}{#2}{#3}}

\def\ddt{{\mathchoice{\FRAC{1}{d}{dt}}%
{\FRAC{1}{d}{dt}}%
{\FRAC{3}{d}{dt}}%
{\FRAC{3}{d}{dt}}}}

\def\ddtm{\ddt^{\tiny -}}

\newcounter{rmnum}
\newenvironment{romannum}{\begin{list}{{\upshape (\roman{rmnum})}}{\usecounter{rmnum}
\setlength{\leftmargin}{14pt}
\setlength{\rightmargin}{8pt}
\setlength{\itemindent}{-1pt}
}}{\end{list}}

\def\itemb{\item[{\small $\bullet$}\ ] }

\def\barr{{\bar r}}

\def\barW{{\overline W}}

\def\barG{{\overline G}}

\def\bfmath#1{{\mathchoice{\mbox{\boldmath$#1$}}%
{\mbox{\boldmath$#1$}}%
{\mbox{\boldmath$\scriptstyle#1$}}%
{\mbox{\boldmath$\scriptscriptstyle#1$}}}}

\def\bfmX{\bfmath{X}}

\def\bfmW{\bfmath{W}}

\begin{document}
\author[krd]{Ken R. Duffy} 
\ead{ken.duffy@nuim.ie}
\author[spm]{Sean P. Meyn}
\ead{meyn@illinois.edu}
\address[krd]{Hamilton Institute, National University of Ireland Maynooth, 
Ireland} 
\address[spm]{
Dept. of Electrical and Computer Engineering,
and the Coordinated Science Laboratory, University of Illinois
at Urbana-Champaign, Urbana, IL, 61801, U.S.A.
}

\title{Most likely paths to error \\
when estimating the mean of  a reflected random walk}

\date{December 22, 2009}

\begin{abstract}
It is known that simulation of the mean position of a Reflected
Random Walk (RRW) $\{W_n\}$ exhibits non-standard behavior, even for
light-tailed increment distributions with negative drift. The
Large Deviation Principle (LDP) holds for deviations below the mean,
but for deviations at the usual speed above the mean the rate
function is null. This paper takes a deeper look at this phenomenon.
Conditional on a large sample mean, a complete sample path LDP
analysis is obtained.  Let $I$ denote the rate function for the one
dimensional increment process. If $I$ is coercive, then given a
large simulated mean position, under general conditions our results
imply that the most likely asymptotic behavior, $\psi^*$, of the paths
$n^{-1} W_{\lfloor tn\rfloor}$ is to be zero apart from on an
interval $[T_0,T_1]\subset[0,1]$ and to satisfy the functional
equation
\begin{align*}
\nabla I\left(\ddt\psi^*(t)\right)=\lambda^*(T_1-t)
\quad
\text{whenever } \psi(t)\neq 0.
\end{align*}
If $I$ is non-coercive, a similar, but slightly more involved,
result holds.

These results prove, in broad generality, that Monte Carlo estimates
of the steady-state mean position of a RRW have a high likelihood
of over-estimation. This has serious implications for the performance
evaluation of queueing systems by simulation techniques where steady
state expected queue-length and waiting time are key performance
metrics. The results show that na\"ive estimates of these quantities
from simulation are highly likely to be conservative.

\end{abstract}



\maketitle

\section{Introduction}
\label{sec:intro}

Consider $\bfmW=\{W_n,n\ge 0\}$, a random walk that starts at zero,
is reflected at the origin, and has increments process
$\bfmX=\{X_n,n\ge0\}$. The Reflected Random Walk (RRW) is governed
by Lindley's recursion \cite{Lindley52},
\begin{align}
\label{eq:lindley}
W_0:=0 \text{ and } W_{n+1} = [W_n+X_n]^+ \text{ for } n\geq1.
\end{align} 
This recursion plays a fundamental r\^ole in queueing systems and
has long been an important object of study in evaluating their
performance. With $\bfmX$ being the difference between service times
and inter-arrival times of customers, the RRW $\bfmW$ describes the
evolution of waiting times at a single server first-come first-served
queue with infinite waiting space. Lindley's recursion also governs
the evolution of the queue-length at certain single server queues,
such as the M/M/1 queue \cite{Asmussen03}.

Since the 1980s, large deviation techniques have been brought to
bear on the analysis of equation \eqref{eq:lindley} and the
distribution of an element of its stationary solution, which exists
whenever $\bfmX$ is stationary \cite{Loynes62}. For example, using
a one-dimensional large deviations approach it has been established in
broad generality that the stationary distribution possesses logarithmic
asymptotics, see
\cite{Glynn94}\cite{Duffield95}\cite{Duffy03}\cite{Duffy04b}\cite{Lelarge08}
and references therein.  This fact is exploited in the theory
effective bandwidths \cite{Kelly96} and in the development of
on-the-fly estimation schemes from observations of the queueing
behavior for the determination of quality of service performance
metrics
\cite{Duffield95A}\cite{Ganesh96}\cite{Gyorfi00}\cite{Paschalidis01}\cite{Duffy05}\cite{Mandjes09}.
Moreover, through the use of functional large deviation techniques,
assuming $\bfmX$ is i.i.d., the seminal paper \cite{Anantharam89}
proved the significant, broadly applicable result that the most
likely path to a large value of the transient RRW is piece-wise
linear. This deduction has since been extended (e.g. \cite{Puhalskii95}),
including the establishment of results in the stationary regime
(e.g. \cite{Dobrushin98}\cite{Ganesh02}\cite{Duffy08}). All of these
papers report piece-wise linear most likely paths to a large value
of the RRW or an element of its stationary solution.

The exclusive focus of all of the research cited above is to garner
understanding of likelihood of large values of the RRW either in
the transient or stationary regime, and the determination of the
most likely paths to these large values. In the present article we
employ a functional large deviation approach to analyze the estimation
of a fundamental quantity for the performance evaluation of a RRW
that has so far been overlooked: its mean value. This study reveals
substantially richer structure than the study of large values of
the RRW, leading to non-convex rate functions and concave most
likely paths. Despite this, perhaps surprisingly, general qualitative
and quantitative deductions can still be made.

\medskip

The starting point of the present paper is the following qualitative result:
It was observed recently  that simulation of the
mean position of a RRW,
\begin{align*}
\barW_n\eqdef\frac{1}{n} \sum_{i=1}^{n} W_i,
\end{align*}
exhibits non-standard behavior, even for light-tailed increments
with negative drift. For example, if $\bfmX$ is i.i.d., then the
probability that $\barW_n$ underestimates the long run expected
value decays exponentially in $n$, but the probability of an
over-estimate decays sub-exponentially. This is shown in the following
proposition, which is taken from \cite{CTCN}. Part (i) follows from
Theorem~11.2.3 and part (ii) from Proposition~11.3.4 (see also
\cite{mey06a}).
\begin{proposition}
\label{thm:MM1LDPfails}
Consider the RRW where $\bfmX$ is i.i.d.\ with $\Expect[X_0]<0$ and
$\Expect[X_0^2]<\infty$. Then the Markov chain $\bfmW$ has a unique
invariant probability measure with finite steady-state mean $\barW$,
and the simulated averages have the following properties.
\begin{romannum}
\item
The lower error-probability decays exponentially:  For each $r<\barW$,
\begin{align*}
\limsup_{n\to\infty}\frac{1}{n}\log  \Prob \bigl\{ \barW_n \le r\bigr\}<0.
\end{align*}
\item 
The upper error-probability decays sub-exponentially:  For each $r>\barW$,
\begin{align*}
\lim_{n\to\infty}\frac{1}{n}\log  \Prob \bigl\{ \barW_n \ge r\bigr\}=0.
\end{align*}
\end{romannum}
\end{proposition}

This paper takes a deeper look at the latter phenomenon, providing
a detailed understanding of why it is hard to simulate the mean
position of a RRW and, therefore, why care must be taken drawing
deductions regarding average queueing performance from the output
of a simulation.  We establish that, in broad generality, the process
$\{n^{-2}\sum_{i=1}^nW_i\}$ satisfies a Large Deviation Principle
(LDP) with a non-trivial rate function.  As a consequence, the
likelihood that the sample-mean estimate of a RRW is an overestimate
decays on a slower than exponential scale.  The rate function in
question is non-convex and this LDP could not, therefore, be
established by asymptotic analysis of scaled cumulant generating
functions, an approach commonly employed in queueing theory and
used in the G\"artner-Ellis method. Unlike the most likely paths
of the RRW that lead to a large position which are piece-wise linear
\cite{Anantharam89}\cite{Ganesh04}, we ascertain that that the most
likely paths associated with a large simulated mean possess more
complex features: they are concave, with a possible discontinuity
when the path first becomes deviant. A number of examples of these
general results are presented to demonstrate the range of qualitative
possibilities.

The results contained in this article clearly indicate that significant
statistical care must be taken when using estimates from simulation
of the mean position of a RRW. This has serious implications for
the performance evaluation of queueing systems by simulation
techniques where steady state expected queue-length and waiting
time are key performance metrics. Our results show, in broad generality,
that the most natural estimation scheme, Monte Carlo estimates, of
these expected values suffer a likelihood of over-estimation that
is, approximately speaking, Weibull-like with shape parameter $1/2$.
Consequently, na\"ive estimation of these quantities from simulation
is likely to underestimate system performance. 

As a concrete illustration of these general results, one example
in which the rate function and most likely paths are explicitly
computable can be found in the following proposition.
\begin{proposition}
\label{thm:likelyGauss}
Consider the RRW in which the increments process $\bfmX$ consists
of i.i.d Gaussian random variables with mean $-\delta<0$ and variance
$\sigma^2$. Then $\{n^{-1}\barW_n\}$ satisfies the LDP in $[0,\infty)$
with rate function
\begin{align}
I_{\barW}(z)
        &= 
        \begin{cases}
	\displaystyle
        \frac{4\delta}{\sigma^2}\sqrt{\frac{z\delta}{6}}
                & \text{if } z\in[0,\delta/6], \\
	\displaystyle
        \frac{3}{2\sigma^2}\left(z+\frac{\delta}{2}\right)^2 
                & \text{if } z\in [\delta/6,\infty). \\
        \end{cases}
\label{eq:Gaussrf}
\end{align}
As $n$ tends to infinity, 
the most likely paths of $n^{-1}W_{\lfloor nt\rfloor}$ leading to
$\barW_n\geq nz$, which we denote $\psi^*$ are as follows.
\begin{romannum}
\item
If $z\in (0,\delta/6]$, then
for any $T_0\in[0,1-\sqrt{6z/\delta}]$
\begin{align*}
\psi^*(t) =
\begin{cases}
0  & \text{ if } t\in[0,T_0]\cup[T_0+\sqrt{6z/\delta},1], \\
\displaystyle
\delta(t-T_0)-\delta \sqrt\frac{\delta}{6z} (t-T_0)^2
	&
	 \text{ for } t\in [T_0,T_0+\sqrt{6z/\delta}].
\end{cases}
\end{align*}
\item
If $z\in [\delta/6,\infty)$, then 
\begin{align*}
\psi^*(t) =
3\left(z+\frac{\delta}{2}\right)\left(t-\frac{t^2}{2}\right) 
	-\delta t \quad \text{ for } t\in [0,1].
\end{align*}
\end{romannum}
\end{proposition}
We return to this example in Section~\ref{sec:examples1} where the
rate function in equation \eqref{eq:Gaussrf} and two most likely
paths are illustrated in Figure~\ref{fig:Gaussrf}.

\Proposition{likelyGauss} and other results that follow concern
asymptotics of the doubly scaled sum  $n^{-1}\barW_n =  n^{-2}
\sum_{i=1}^{n} W_i$. The $n^2$ scaling is similar to \cite[Theorem
4.1]{Borovkov03}, concerning asymptotics for the GI/G/1 queue in
the light tailed setting.  This result proves that the tail of the
busy time distribution decays more slowly than exponentially:
$\lim_{n\to\infty}n^{-1}\log P\left(B> n^2z\right)=-K\sqrt{z}$ for
each $z>0$ and some $K>0$, where $B$ denotes the busy time in
steady-state.   The form of the limit can be predicted through
scaling arguments: If $\{n^{-2}B\}$ satisfies the LDP, it must do
so with a rate function of the form $K\sqrt{z}$, as can be seen by
considering the substitution $m=n\sqrt{z/y}$. The rate function for
the large deviations of the process $\{n^{-1}\barW_n\}$ is necessarily
more complex and, in general, the rate function will diverge more
rapidly than $\sqrt{z}$ as $z\to\infty$.

The rest of this paper is organized as follows.
In Section~\ref{sec:ldp} we prove in broad generality that the
sample paths of the rescaled simulated mean $n^{-1}\barW_n$ satisfy
a functional LDP. Using this LDP, in Section~\ref{sec:mostlikelypaths}
we characterize properties of the most likely paths of the RRW given
that the rescaled simulated mean is large. In Section~\ref{sec:examples}
we present examples including the RRW with i.i.d. Gaussian increments,
the M/D/1 queue, the D/M/1 queue and the M/M/1 queue, as these
exhibit the full range of theoretically possible behaviors.

\section{Functional LDP for the rescaled simulated mean position}
\label{sec:ldp}

We assume that the reader is familiar with the basics of large
deviation theory, such as the definition of the LDP and the statement
of the Contraction Principle, as can be found in
\cite{Deuschel89}\cite{Dembo98}\cite{denHollander00}\cite{Ganesh04}. The
notation in this paper is as follows. Let $C[0,1]$ denote the set
of continuous $\R$-valued functions on $[0,1]$ equipped with the
topology induced by the supremum norm,
$\|\phi\|=\sup_{t\in[0,1]}|\phi(t)|$. Let $D[0,1]$ denote the space
of $\R$-valued c\'adl\'ag functions on $[0,1]$ equipped with the
Skorohod ($J_1$) topology
\cite{Skorohod56}\cite{Billingsley68}\cite{Whitt02} induced by the
following metric: for any two functions $\phi,\psi\in D[0,1]$,
define
\begin{align*} 
d(\phi,\psi) := \inf_{\lambda\in\Lambda} 
\bigl\{
        \max\left(\|\phi\circ\lambda-\psi\|,\|\lambda-e\|\right)
        \bigr\},
\end{align*} 
where $e$ is the identity ($e(t) = t$),   and $\Lambda$ is the set of
strictly increasing functions $\lambda$ from $[0,1]$ to $[0,1]$
that are continuous, with a continuous inverse.
Finally, let $L[0,1]\subset
D[0,1]$ denote the set of functions that have finite variation.
Each $\phi\in L[0,1]$ has a Lebesgue decomposition with respect to
Lebesgue measure whose absolutely continuous part we denote
$\phi^{(a)}$ and whose singular component we denote $\phi^{(s)}$,
so that $\phi(t) = \int_0^t \dot\phi^{(a)}(s) ds +\phi^{(s)}(t)$.
Furthermore, we decompose $\phi^{(s)}$ into its positive $\phisu$
and negative $\phisd$ parts by the Hahn Decomposition Theorem.

For each $n\in\N$ and all $t\in[0,1]$, we define the following
scaled sample paths:
\begin{align*}
x^n(t) := \frac 1n \sum_{i=0}^{\lfloor nt \rfloor -1} X_i,
\ \ 
w^n(t) := \frac 1n W_{\lfloor nt \rfloor}  
\ \ 
\text{and}
\ \ 
\wbar^n(t) := 
	\frac{1}{n^2}\sum_{i=1}^{\lfloor nt\rfloor} W_i
		+\frac{1}{n^2}(nt-\lfloor nt\rfloor])W_{\lfloor nt\rfloor+1}.
\end{align*}
The first two of these are elements of $D[0,1]$ and correspond to
the paths for the simulated position of the unconstrained random
walk and for the simulated position of the RRW, respectively.  The
sample path $\wbar^n$ is an element of $C[0,1]$ and is the polygonally
approximated continuous path for the rescaled simulated mean location
of the RRW. In particular, note that
$\wbar^n(1)=n^{-1}\barW_n=n^{-2}\sum_{i=1}^n W_i$ is the rescaled
simulated mean of the RRW.

For the general qualitative theorem we make the following assumption.
\begin{assumption}
\label{ass:LDPX}
The sample paths for the unconstrained random walks $\{x^n\}$
satisfy the LDP in D[0,1] with good rate function $I_X$.
\end{assumption}

This assumption is known to hold for a large collection of processes.
For example, if $\bfmX$ is an i.i.d. sequence, then define $\thetad:=
\sup\{\theta>0: \Expect[\exp(-\theta X_0)]<\infty\}$ and $\thetau:=
\sup\{\theta>0: \Expect[\exp(\theta X_0)]<\infty\}$.  If
$\min\{\thetad,\thetau\}>0$, then by Cram\'er's Theorem
\cite{Cramer38}\cite{Dembo98} the partial sums of $\{x^n(1)\}$ satisfy
the LDP in $\R$ with the good, convex (local) rate function
\begin{align}
\label{def:Il}
I(y) := \sup_{\theta}\left(\theta y -\log \Expect[\exp(\theta X_0)]\right),
\text{ for } y\in\R, 
\end{align}
and Mogul'skii's Theorem \cite{Mogulskii76} proves that Assumption
\ref{ass:LDPX} holds true. The rate function is typically of the form
(e.g 
\cite{Lynch87}\cite{Mogulskii93}\cite{Puhalskii95}\cite{Puhalskii98}):
\begin{align}
\label{eq:IDZ}
I_X(\gamma) = \begin{cases}
        \int_{0}^1 I(\dot \gamma^{(a)}(s)) ds
        + \thetad\gammasd(1) + \thetau\gammasu(1)
                &\; \mbox{if}\;
                        \gamma\in L[0,1],\\
        +\infty &\; \mbox{otherwise},
        \end{cases}
\end{align}
where $\thetad\gammasd(1):=0$ if $\thetad=\infty$ and $\gammasd(1)=0$, and
$\thetau\gammasu(1):=0$ if $\thetau=\infty$ and $\gammasu(1)=0$.  Dembo and
Zajic \cite{Dembo95} have generalized Mogul'skii's Theorem to include
sequences $\bfmX$ that need not be i.i.d, but that satisfy a uniform
super-exponential tail condition that ensures that the generalization
of $\min(\thetad,\thetau)$ is $+\infty$, as well as a mixing condition
that encompasses, for example, Markov chains that are uniformly
ergodic. The resulting rate function for these processes is also
of the form in equation \eqref{eq:IDZ}, but the cumulant generating
function $\log \Expect[\exp(\theta X_0)]$ in equation \eqref{def:Il}
is replaced with the scaled cumulant generating function  $\lim
n^{-1}\log \Expect[\exp(\theta(X_0+\cdots+X_{n-1}))]$.

\begin{theorem}
\label{thm:waiting}
The following hold under Assumption \ref{ass:LDPX}:
\begin{romannum}
\item
The sequence of rescaled
paths of the simulated mean of the RRW $\{\wbar^n\}$ satisfies the
LDP in $C[0,1]$ with rate function
\begin{align*}
I_{\barW}(\phi)=\inf_{\gamma \in L[0,1]}\left\{I_X(\gamma):
\int_0^t \sup_{s\leq t}\left(\gamma(t)-\gamma(s)\right)\,dt 
	= \phi(t)\; \text{ for all }\; t\in[0,1]\right\}.
\end{align*}
\item
If $I_X$ is of the form in equation \eqref{eq:IDZ}, then $I_{\barW}$
is only finite at those functions $\phi$ such that $\dot\phi$ exists, $\dot\phi$
is non-negative and $\dot\phi$ is an element of $L[0,1]$, in which
case
\begin{align}
	\label{eq:IWintD}
I_{\barW}(\phi)=
	\int_0^1 
	\left(
	I(\ddot \phi^{(a)}(s)) 1_{\{\dot \phi(s)>0\}}
	+ \inf_{y\leq0}I(y) 1_{\{\dot \phi(s)=0\}}
	\right)\,ds
		+\thetad\phidsd(1)
		+\thetau\phidsu(1).
\end{align}
\end{romannum}
\end{theorem}
\pf{Proof}
The proof of the first assertion follows from the contraction
principle (e.g. \cite[Theorem 4.2.16]{Dembo98}) after noting the
following. The Skorohod map, $f(\gamma)(t)=\gamma(t)-\inf_{s\leq
t} \gamma(s)$, is continuous from $D[0,1]$ to $D[0,1]$ (e.g.
\cite[Theorem 13.5.1]{Whitt02}) and $f(x^n)(t) = w^n(t)$. The
integration map, $g(\psi)(t)= \int_0^t\psi(s)\, ds$ is continuous
from $D[0,1]$ to $C[0,1]$ (e.g. \cite[Theorem 11.5.1]{Whitt02}) and
$g(w^n)(t)=\wbar^n(t)$.

For the second assertion, if $\phi$ is such that $\dot\phi$ does
not exist, takes negative values or is not an element of $L[0,1]$,
then $I_\barW(\phi)=+\infty$, as can be seen from the contraction
principle. If $\dot\phi$ exists, is non-negative and is an element
of $L[0,1]$, then
\begin{align*}
I_{\barW}(\phi)
        &= \inf_{\gamma\in L[0,1]}
                \left\{\int_0^1 I(\dot \gamma^{(a)}(s))\, ds
                +\thetad \gammasd(1)
                +\thetau \gammasu(1)
                : \right. \\
        &\qquad\qquad\qquad\left.
        \sup_{s\leq t}\left(\gamma(t)-\gamma(s)\right)
        = \dot\phi(t)\; \mbox{for all}\; t\in[0,1]\right\}.
\end{align*}
If $\dot\phi(t)>0$, then $\gamma$ must satisfy
$\ddot\phi^{(a)}(t)=\dot\gamma^{(a)}(t)$. As $\ddot\phi^{(a)}(t)=0$
for almost all $t$ such that $\dot\phi(t)=0$, if $\dot\phi(t)=0$
we are free to choose $\dot\gamma^{(a)}(t)=\inf_{y\leq0}I(y)$ to
minimize the rate function. The singular parts $\phidsd$ and $\phidsu$
must be mimicked by $\gamma^{(s)}$. In order to minimize the rate
function, $\gamma^{(s)}$ is unchanging everywhere else, leading to
the result.
\qed 

The rate function in equation \eqref{eq:IWintD} can be understood
as follows. In order to see the rescaled simulated mean sample path
$\phi$, in the integral one must locally pay for changes in the
increments process so long as the location is positive. If the
location is zero, then the increments can take their most likely
value less than or equal to zero. The singular parts of the location
are matched by singular parts in the increments.

\section{Most likely RRW paths to a large simulated mean position}
\label{sec:mostlikelypaths}

Considering Theorem \ref{thm:waiting} in conjunction with the
contraction principle and the projection $\phi\mapsto\phi(1)$,
roughly speaking, we can deduce that
\begin{align*}
\Prob\left\{\frac{1}{n^2}\sum_{i=1}^nW_i \approx z\right\} 
	\sim 
	\exp\left(-n \inf_{\phi\in C[0,1]}
	\left\{I_\barW(\phi):\phi(1)=z\right\}\right).
\end{align*}
Thus consider the following minimization problem:
\begin{align*}
\inf_{\phi\in C[0,1]} \left\{I_\barW(\phi):\phi(1)=z\right\}.
\end{align*}
If $I_X$ is of the form in equation \eqref{eq:IDZ}, then this problem
can be rewritten in terms of the fluid limit paths of the RRW:
\begin{equation}  
\begin{aligned}
\text{{\bf minimize} } \quad& J(\psi) \\
\text{{\bf subject to} } \quad& 
\psi\in L^+[0,1] \text{ and } \int_0^1 \psi(s)\, ds = z,
\end{aligned}
\label{eq:main}
\end{equation}
where $L^+[0,1]$ is the set of non-negative elements of $L[0,1]$
and the objective function is
\begin{align}
\label{def:J}
J(\psi) \eqdef\int_0^1 
	\left(
	I(\dot\psi^{(a)}(s)) 1_{\{\psi(s)>0\}}
	+ \inf_{y\leq0}I(y) 1_{\{\psi(s)=0\}}
	\right)\,ds
		+\thetad\psisd(1)
		+\thetau\psisu(1).
\end{align}
The evaluation of the optimization \eqref{eq:main} and the
identification of properties of its infimal argument (or arguments)
are the subject of the rest of this paper. That is,
we wish to identify properties of the most likely fluid simulated
RRW paths that give rise to the simulated mean position being
unusually large. These optimizers are most likely paths in the sense
that if $G$ is any measurable neighborhood of the set of minimizing
arguments to \eqref{eq:main} and $\barG=\{\phi:\phi(t)=\int_0^t\psi(s)\,ds
\text{ for some } \psi\in G\}$, then $\lim_{n\to\infty}P\{\wbar^n\in
\barG\}=1$. This follows, for example, by \cite[Theorem 2.2]{Lewis95}.

In addition to Assumption \ref{ass:LDPX}, the following assumption
is in force throughout the rest of this article.
\begin{assumption}
\label{ass:LDPexplicit}
The rate function $I_X$ is of the form in equation \eqref{eq:IDZ}, where
$I$ is a good, convex rate function, and there exists $\delta>0$
such that $I(-\delta)=0$, so that the RRW is stable.
\end{assumption}

Note that as $I$ is a rate function, it is lower semi-continuous.
The maximal value for which it is finite is denoted by
\begin{equation}
\barr \eqdef \sup\{r : I(r)<\infty\} .
\label{e:barr}
\end{equation}
Suppose that $I$ is a non-coercive function: $\barr<\infty$ and
$\lim_{r\uparrow \barr} I(r)<\infty$.
Then the limit must coincide with $I(\barr)$, which is thus finite.
This is needed to ensure the existence of optimal paths. Note
also that $I$ being non-coercive is mutually exclusive with $\thetau<\infty$,
which requires $I(r)<\infty$ for all positive $r$.

\begin{theorem}
\label{thm:props} 
An optimal solution to the optimization problem \eqref{eq:main}
exists, and any optimal solution $\psi^*$  satisfies the following
properties: There exists $0\le T_0 < T_1\le 1$ such that, 
\notes{K.D. Have altered the language, but not the meaning, at the
start of the proof.}
\begin{romannum}
\item $\psi^*(t)>0 $ on the open interval $(T_0,T_1)$ and
$\psi^*(t)=0$ for $t\in [0,1]\setminus[T_0,T_1]$;
\item
$\psi^*$ is concave on $[T_0,T_1]$;
\item
$\psi^*$ is continuous on $(T_0,T_1]$, with a possible jump at $t=T_0$.
\end{romannum}
\end{theorem}

The proofs of this theorem and the two that follow are postponed to the end of this section.

The time $T_0$ is taken to be the minimal time that a path is
non-zero and $T_1$ the maximal time that it is non-zero:
\begin{equation}
T_0\eqdef\inf\{t\ge 0: \psi(t) >0\}
\ \ 
\text{and}
\ \ 
T_1\eqdef\sup\{t\le 1 : \psi(t) >0\}.
\label{e:T01}
\end{equation}
If $0<\barr<\infty$, then we define 
\begin{align}
T_0^0\eqdef \sup\{t\ge 0: \ddtm\psi(t) =\barr\}.
\label{e:T00}
\end{align}
If the supremum is over an empty set then we take  $T_0^0 =T_0$; hence
the inclusion $T_0^0\in [T_0,T_1]$ follows from the definitions.
The following theorem identifies the structure of the most likely
path for $t$ between $T_0^0$ and $T_1$. The one that follows it identifies
how most likely paths must end.

\begin{theorem}
\label{thm:Calculus} 
Let $0\le T_0\le T_0^0 <T_1\le 1$ denote the values given in
\eqref{e:T01} and \eqref{e:T00}. Then for $\psi^*$ to be an optimal
path, there must exist constants $b\in\Re$ and $\lambda^*>0$ such that
\begin{align}
\nabla I(\dot \psi^*(t)) =b-\lambda^* t  \ \  \text{ for all }  T_0^0<t<T_1.
\label{e:Calc}
\end{align}
In particular, if $I$ is coercive, then $T_0^0=T_0$ and
equation \eqref{e:Calc} is satisfied for all $t$ such that $\psi^*(t)>0$.
\end{theorem}

\begin{theorem}
\label{thm:FinalDerValue}
Suppose that  $T_0^0<T_1 $. Then  $\ddtm\psi^*(t)|_{t=T_1} =-\delta$
and $b=\lambda^*T_1$ in equation \eqref{e:Calc}.
\notes{K.D. Added latest observation to Theorem \ref {thm:FinalDerValue}
plus one supporting line at end of corresponding proof. Statement
of Proposition \ref{thm:reduced} also changed to reflect this.}
\end{theorem}

As well as providing insight into the structure of the most likely
paths, Theorems~\ref{thm:props}, \ref{thm:Calculus} and
\ref{thm:FinalDerValue} enable the reduction of the problem
\eqref{eq:main} from an infinite dimensional optimization to a
finite dimensional optimization problem that can be readily solved
numerically, if not analytically.

\begin{proposition}
\label{thm:reduced}
Given $z>0$, define the subset $S_z\subset L^+[0,1]$ of 
{\em  potential solutions} 
to be the collection of functions $\psid$ such that, for
some $0\le T_0 <T_1\le 1$ and $T_0^0\in [T_0,T_1]$:
\begin{romannum}
\item $\psid(t) = 0$ for $t\in[0,T_0)\cup(T_1,1]$;
\item if $T_1<1$, then $\psid(T_1)=0$;
\item $\int_{T_0}^{T_1}\psid(t)dt =z$.
\item if $T_0^0<T_1$, then $\ddtm\psi^*(t)|_{t=T_1}=-\delta$.
\item $\nabla I(\dot \psi^*(t))=\lambda^* (T_1-t)$ for all
$t\in(T_0^0,T_1)$.
\end{romannum}
If $I$ is coercive and $\thetau=\infty$, then in addition to (i)-(v):
\begin{romannum} 
\itemb  
$T_0^0=T_0$; 
\itemb  
$\psid$ has no discontinuities.
\end{romannum}
If $I$ is coercive and $\thetau<\infty$, then in addition to (i)-(v):
\begin{romannum} 
\itemb 
$T_0^0=T_0$; 
\itemb 
if $\psid$ has a discontinuity, it is at $T_0$.
\end{romannum} 
If $I$ is non-coercive (which ensures that $\thetau=\infty$), then in
addition to (i)-(v):  
\begin{romannum}
\itemb  
$\dot \psid(t) = \barr $ for $t\in[T_0,T_0^0)$;
\itemb  
$\psid$ has no discontinuities.
\end{romannum}
The problem \eqref{eq:main} is then equivalent to
\begin{equation} 
\begin{aligned} 
\text{{\bf minimize} } & \quad
	I(\barr)(T_0^0-T_0)
	+\int_{T_0^0}^{T_1} I(\dot\psid(s)) dt
	+\psid^\uparrow(T_0)\thetau
\\
 \text{{\bf subject to} } & \quad
\psid\in S_z.  
\end{aligned}
	\label{eq:reduced}
\end{equation}
\end{proposition}
After proving Theorems~\ref{thm:props}, \ref{thm:Calculus} and
\ref{thm:FinalDerValue}, in Section~\ref{sec:examples} we will use
the reduced representation of the problem defined in equation
\eqref{eq:reduced} in the consideration of illustrative examples.

The following lemma will be used to establish properties of an
optimal fluid trajectory.
\begin{lemma}
\label{thm:psiLower}
Suppose that $\psi^0$ is a fluid trajectory satisfying  
$z(0)\eqdef \int_0^1 \psi^0(s)\, ds <\infty$.
\begin{romannum}
\item
For any $d\ge 0$ and $t\in[0,1]$ define $\psi^d(t) =\max(0,
\psi^0(t)-d)$, and $z(d)\eqdef \int_0^1 \psi^d(s)\, ds$. Then
$z(\varble)$  is convex and non-increasing as a function of $d$.
\item
If $\psid^\uparrow(1)=0$ then $J(\psi^d)$ is non-increasing as a function of $d$. 
\end{romannum}
\end{lemma}

\pf{Proof}
For each $t$, the function of $d$ given by $\psi^d(t)=\max(0,\psi^0(t)-d)$
is convex and non-increasing. It follows that its integral over
time is also concave and non-increasing.

Part (ii) then follows from the definition of $J$ given in equation
\eqref{def:J}.
\qed

\pf{Proof of \Theorem{props}.}
We first establish the existence of an optimizer, which follows
from topological arguments. The objective function $J:D[0,1]\mapsto[0,1]$
is defined for elements of $L^+[0,1]$ in equation \eqref{def:J};
set $J(\psi) = +\infty $ for $\psi\not\in L^+[0,1]$. The function
$J$ is lower semi-continuous and has compact level sets as it is
the good rate function for the LDP of the sample path process
$\{w^n\}$. With domain $ D[0,1]$, the mapping
$\psi\mapsto\int_0^1\psi(s)\,ds$ is continuous
(e.g. \cite[Theorem 11.5.1]{Whitt02}), so that the set $\{\psi\in
D[0,1]:\int_0^1\psi(s)\,ds =z\}$ is closed. In a Hausdorff space,
the infimum of a lower semi-continuous function with compact level
sets is attained on closed sets (e.g. \cite[Lemma~4.1]{Ganesh04}).
Thus if the infimum in \eqref{eq:main} is finite, it is attained
at some $\psi^*\in L^+[0,1]$ such that $\int_0^1\psi^*(s)\,ds=z$.

Regarding the properties of an optimizer $\psi^*$, note first that
it is obvious that $\psisd(1)=0$: By removing downward jumps we
reduce $J(\psi)$, while increasing the area $\int_0^1 \psi(s)\,
ds$. On letting $\psi^0$ denote the new trajectory, and setting
$z(d)\eqdef \int_0^1 \psi^d(s)\, ds$, \Lemma{psiLower} then implies
that $z(d)=z$ for some $d\ge 0$, with $J( \psi^d)\le J(\psi)$.

We assume without loss of generality that the closure of $\{t : \psi(t)>0\}$ is equal to the interval $[T_0,T_1]$ (where the endpoints are defined in \eqref{e:T01}):     If there exist times $t_0<t_1$ satisfying $t_0\ge T_0$,  $t_1\le T_1$, and $\psi(t)=0$ for $t\in (t_0,t_1)$, then the trajectory can be shifted as follows,
\[
\psi^0(t) = \begin{cases} \psi(t) & t\in [T_0,t_0]
				\\
				\psi(t+(t_1-t_0))  & t\in [t_0,1-   (t_1-t_0)]
				\\
				\max(0,\psi(1) + (1-   (t_1-t_0) - t) \delta)    & t\in (1-   (t_1-t_0) ,1).
		\end{cases}
\]
Once again, on setting $z(d)\eqdef \int_0^1 \psi^d(s)\, ds$, we have $z(0)\ge z$,  and on applying \Lemma{psiLower} we have $z(d)=z$ for some $d\ge 0$, with $J( \psi^d)\le J(\psi)$. 

Next,  we assume without loss of generality that $T_0=0$:  We can replace $\psi$ by,
\begin{align*}
\psi^0(t) = \begin{cases} \psi(t-T_0) & t\in [0,1-T_0]
\\
				\max(0,\psi(1-T_0)  +  (1-T_0-t) \delta )    & t\in (1-T_0,1).
		\end{cases}
\end{align*}
An application of   \Lemma{psiLower} again shows that $J(\psi^d) \le J(\psi^0) =  J(\psi)$, and $z(d)=z = \int_0^1 \psi(t)\, dt$ for some $d\ge 0$. 

\begin{figure}[h]
\Ebox{1}{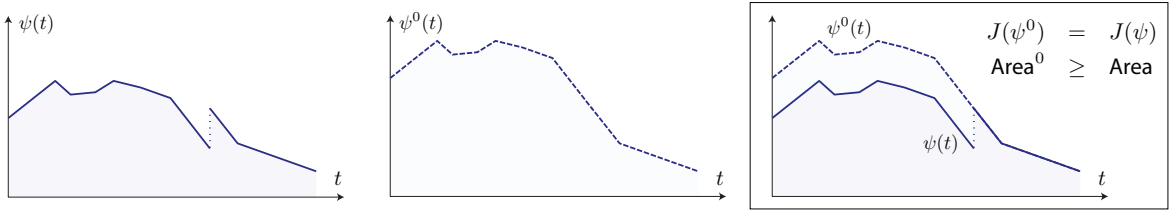}
\caption{A jump for $t>T_0$ cannot be optimal. The rate functional
evaluated at the two paths  $\psi$ and $\psi^0$ is the same,
yet the area is greater using $\psi^0$.} 
\label{fig:JumpPL}
\end{figure}

We can now prove (iii):  \Figure{JumpPL} illustrates why a jump following time $T_0$ cannot be optimal.   A formal proof can be performed through construction as in the previous steps. We define, for any $\psi$,  the new trajectory with $\psi^0(0)=\psi(0)$, and
\[
\psi^0(t) = \psi^\uparrow(1) + \psi(t),\qquad 0<t\le 1.
\]
We have $J(\psi^0) =J(\psi)$, and $ \int_0^1 \psi^0(s)\, ds  \ge  \int_0^1 \psi(s)\, ds$.  Applying \Lemma{psiLower} we have $z(d)=z$ for some $d\ge 0$, with $J( \psi^d)\le J(\psi)$. This proves (iii).

\begin{figure}[h]
\Ebox{1}{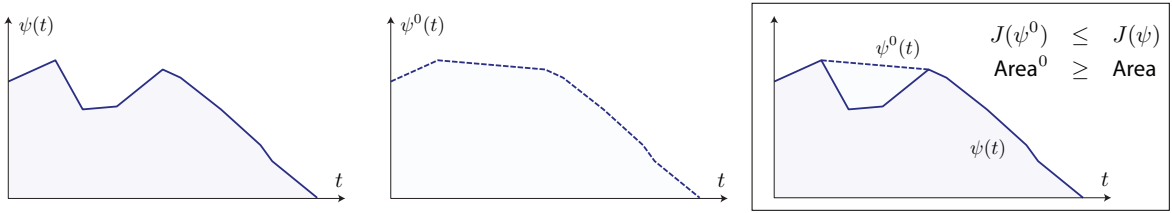}
\caption{An optimal path is concave on $(T_0,T_1)$.   Convexity of $I(r)$ implies that $J(\psi^0)\le J(\psi)$, yet the area is greater using $\psi^0$. }
\label{fig:ConcavePL}
\end{figure}

Similar reasoning establishes concavity of an optimal path --- \Figure{ConcavePL} shows a transformation of a given trajectory to form a new trajectory with reduced value $J(\psi^0)$, but strictly greater area.  Applying \Lemma{psiLower} we obtain (ii).   Part (i) then follows from (ii).
\qed 

\pf{Proof of \Theorem{Calculus}.}
From Theorem \ref{thm:props}, if $\psi^*$ is an optimal solution for \eqref{eq:main} then
it is continuous apart from at $T_0$; if $T_1<1$, then $\psi^*(T_1)=0$.
Thus \eqref{eq:main} can be considered as identifying
\begin{align}
\label{eq:prob1}
&
\inf_{\psi(T_0)\geq0}\inf_{T_1\leq 1}\inf_{T_0^0\le z/\barr}
\left\{(T_0^0-T_0)I(\barr)+\int_{T_0^0}^{T_1} I(\dot\psi(t)) dt
	+\thetau \psi(T_0)
	\right. \\ &\qquad\qquad\qquad: \left.
	\, (T_0^0-T_0)^2\frac{\barr}{2}+ (T_1-T_0)\psi(T_0)
		+\int_{T_0^0}^{T_1}\psi(t)dt = z
		,\,  \psi(T_1)=0 \text{ if } T_1<1 \right\},\nonumber
\end{align}
For fixed $T_0$, $\psi(T_0)$, $T_1$ and $T_0^0$, we are left to
consider finding the solution of a problem of the following kind:
\[\begin{aligned} 
\text{{\bf minimize} } & \quad
	\int_0^T I(\dot\psi(t))\, dt\\
 \text{{\bf subject to} } & \quad \int_0^T \psi(s)\, ds = z'.
\end{aligned}
\]
If $\psi$ is feasible path, then integration by parts gives
\begin{equation}
\int_0^T t \dot \psi(t)\, dt =T \psi(T) -z'. 
\label{e:psiConstraint}
\end{equation}

Introduce the Lagrangian
\begin{align*}
\clL(\psi,\lambda) = \int_0^T I(\dot \psi(t)) \, dt  
	+  \lambda \Bigl(\int_0^T t \dot \psi(t)\, dt+z'
		-T\psi(T)\Bigr).
\end{align*}
There exists $\lambda=\lambda^*$  so that complementary slackness
holds. Hence the optimizer  $\psi^*$ of \eqref{eq:main} also minimizes
$ \clL(\psi,\lambda^*)$ over all $\psi$.  The constant $  \lambda^*$
exists by \cite[(Theorem 1 of Section 8.3)]{lue97},  which only
requires feasibility of  \eqref{e:psiConstraint}  for $z''$ in a
neighborhood of $z'$ (which is true when $T_0^0<T_1$. If $T_0^0=T_1$
then there is nothing to prove).

If $\psi^*$ minimizes the Lagrangian, and if $\delta$ represents a
perturbation satisfying $\delta(t)=0$ for $t\in (T_0^0,T_1)^c$,
then,
\begin{align*}
0=
\frac{d}{d\theta}\clL(\psi^*+\theta \delta,\lambda^*) \Big|_{\theta=0}
=
\int_0^T[ \nabla I(\dot \psi(t)) + \lambda^* t ] \dot \delta (t)\, dt  
\end{align*}
It follows that there exists a constant $b$ such that 
\begin{align*}
\nabla I(\dot \psi(t)) =b-\lambda^* t 
\quad
\text{ for {\it a.e.} $t\in (T_0^0,T_1)$. }
\end{align*}
Returning to problem \eqref{eq:prob1}, this implies that irrespective
of the optimal values of $T_0$, $\psi(T_0)$, $T_0^0$ or $T_1$, the
optimal path satisfies $\nabla I(\dot\psi^*(t))=b-\lambda^*t$ for
$t \in (T_0^0,T_1)$.
\qed 

\begin{figure}[h]
\Ebox{1}{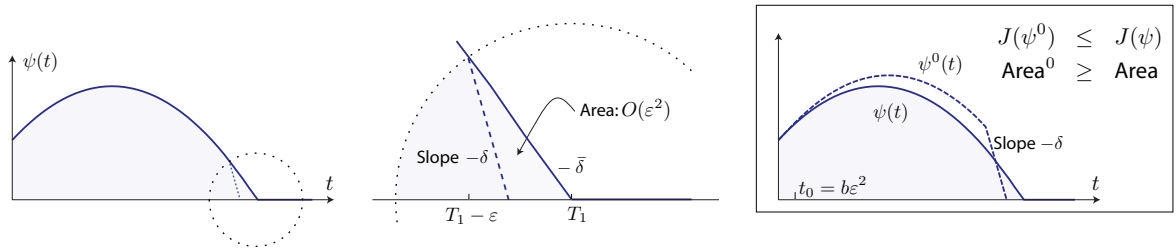}
	\vspace{-.75cm}
\caption{An unnatural slope is costly near $t=T_1$.}
\label{fig:NaturalEndingLong}
\end{figure}

\pf{Proof of \Theorem{FinalDerValue}.}
\Figure{NaturalEndingLong} illustrates the idea of the proof in the special case $T_1<1$.
Let $\bar\delta = - \ddtm\psi(t)|_{t=T_1}$, and suppose that $\bar\delta<\delta$.  We will construct a new trajectory $\psi^0$ with a lower value of $J(\psi^0)$, and increased area.  An application of \Lemma{psiLower} will then show that $\psi$ cannot be optimal.

We henceforth assume without loss of generality that $T_0=0$. 
Concavity of $\psi$ on $[0,T_1]$ implies that   $\ddt \psi(t) \ge -\bar\delta$ for all $t\in (T_0,T_1)$.      
The main idea of the proof (as illustrated in the figure) is as follows: For given $\epsy>0$, the cost contribution to $J(\psi)$ over $[T_1-\epsy,T_1]$ is greater than $\epsy I(-\bar\delta)=O(\epsy)$ if $\bar\delta<\delta$.    However, the additional area obtained is bounded by  $O(\epsy^2)$.

We let $-\delta^0 $ denote the right derivative, following a possible jump:
\[
\delta^0=- \lim_{t\downarrow 0}\ddt\psi(t)
\]
For fixed $\epsy>0$, $b>0$,  let $t_0=b\epsy^2$, and let $\psi^0$ denote the concave function defined by $\psi^0(0+)=\psi(0+)$,  with derivatives for $t>0$ defined as follows:
\begin{align}
\ddt\psi^0(t) = \begin{cases} -\delta^0 & t\in (0,t_0]
				\\
				\ddt\psi(t-t_0)  & t\in (t_0, T_1+t_0 -\epsy]
				\\
				-\delta   & t>T_1+t_0 -\epsy  \ \ \hbox{(provided $\psi^0(t)>0$).}	\end{cases}
\end{align}
For  $b$ sufficiently large, we have $\int_0^1 \psi^0(s)\, ds \ge  \int_0^1 \psi(s)\, ds $ for all $\epsy>0$ sufficiently small.   For the same constant $b$ we also have $J(\psi^0) \le J(\psi) -O(\epsy) + O(\epsy^2)$, so that $J(\psi^0) \le J(\psi)$ for sufficiently small $\epsy>0$.   Fixing $b$ and $\epsy$ so that these bounds hold,  we then apply
\Lemma{psiLower} to conclude that   $z(d)=z$ for some $d\ge 0$, with $J( \psi^d)\le J(\psi^0) \le J(\psi)$.  
The second statement of the theorem follows from $\ddtm\psi^*(t)|_{t=T_1}=-\delta$
and Theorem \ref{thm:Calculus} on noting that $\nabla I(-\delta)=0$.

\qed

\section{Examples}
\label{sec:examples}

\subsection{Coercive rate function: continuous paths}
\label{sec:examples1}

We present two examples with coercive rate functions. One is the
RRW with Gaussian increments. Here the rate function and most likely
path can be determined explicitly. The RRW in the second example
corresponds to the queue-length at departures of an M/D/1 queue
with batch services. In this case, identification of the rate
function and the most likely paths requires the solution of two
transcendental equations, which can be readily obtained numerically.

{\bf Gaussian increments}.
Let $\bfmX$ be i.i.d. Gaussian, with $X_0$ having mean $-\delta<0$
and variance $\sigma^2$. Then the conditions of Mogul'skii's Theorem
are met with $\thetad=\thetau=+\infty$ and the local rate function is
\begin{align*}
I(x) = \frac{1}{2\sigma^2} (x+\delta)^2.
\end{align*}
As $\thetad=\thetau=+\infty$ the sample path rate function $I_\bfmX$
is only finite at absolutely continuous functions.

Without loss of generality, assume that $T_0=0$ and define $T=T_1$
(so that $T$ represents $T_1-T_0$). By Proposition \ref{thm:reduced},
to solve the problem \eqref{eq:reduced} for a given $z>0$, for each
$T\in(0,1]$ we first identify paths $\psid$ that satisfy $\nabla
I(\dot\psid(t))=\lambda^*(T-t)$ in $[0,T]$ and $\dot\psid(T)=-\delta$,
which leads to candidate solutions satisfying
\begin{align*}
\dot\psid(t) = \sigma^2\lambda^*(T-t)-\delta.
\end{align*}
If $T<1$, then in addition we have that $\psid(T)=0$ and
$\int_0^T\psid(t)\,dt = z$ giving $\sigma^2\lambda^* = 2\delta/T$ and
$T=\sqrt{6z/\delta}$. Note that $T<1$ only if $6z<\delta$ and
therefore the optimal path is 
\begin{align}
\label{eq:examplegauss}
\psi^*(t) = \delta t - \delta\sqrt{\frac{\delta}{6z}}\frac{t^2}{2}
\text{ for } t\in[0,T] \text{ if } z<\delta/6.
\end{align}
If $T=1$, then  we have that $c\eqdef \psid(1)\geq0$ and
$\int_0^1\psid(t)\,dt=z$, giving $\sigma^2\lambda^* =2(c+\delta)$ and
$c=3/2(z-\delta/6)$. Note that $c\geq0$ only if $6z\geq \delta$ and
therefore the optimal path is
\begin{align}
\label{eq:examplegauss1}
\psi^*(t) = 
	3\left(z+\frac{\delta}{2}\right)\left(t-\frac{t^2}{2}\right)
		-\delta t
\text{ for } t\in[0,1] \text{ if } z\geq \delta/6.
\end{align}
Evaluating $\int_0^TI(\dot\psi^*(t))\,dt$, for $\psi^*(t)$ defined
in equations \eqref{eq:examplegauss} and \eqref{eq:examplegauss1},
we obtain the rate function $I_\barW(z)$ presented in equation
\eqref{eq:Gaussrf}, Proposition \ref{thm:likelyGauss}, which can be
found in Section~\ref{sec:intro}.

\begin{figure}[ht]
\EboxRatePath{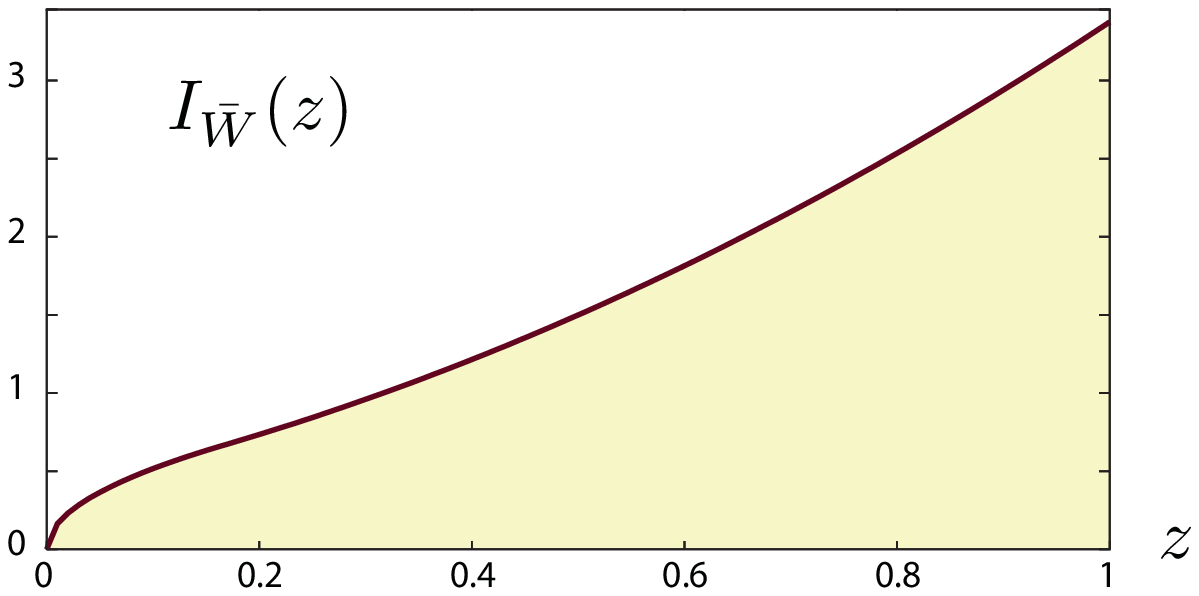}{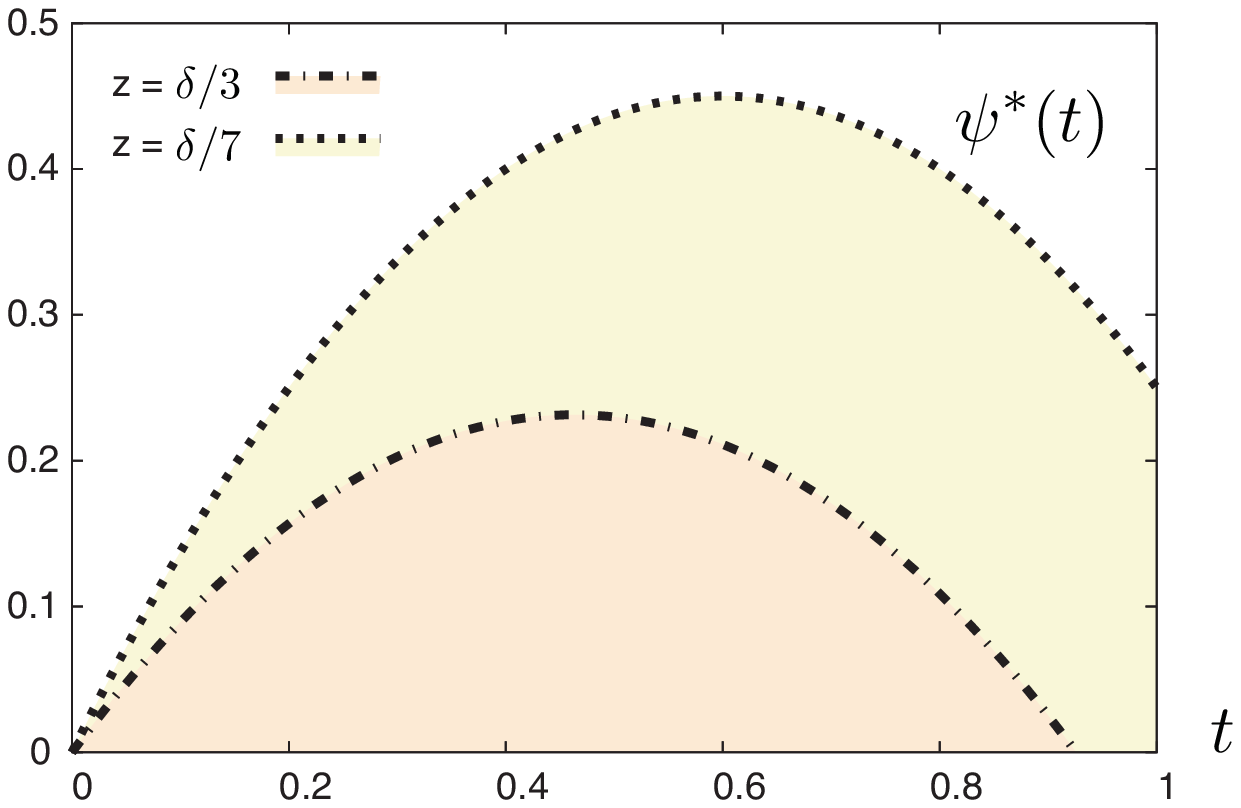}
\caption{Coercive rate function example: i.i.d. Gaussian increments.
On the left hand side is $I_\barW(z)$ versus $z$. Shown on
the right hand side are most likely RRW paths,
$\psi^*(t)\approx n^{-1}W_{\lfloor nt \rfloor}$,
given that $\barW_n\approx nz$.}
\label{fig:Gaussrf}
\end{figure}

With $\delta=1$ and $\sigma^2=1$, the rate function defined in
equation \eqref{eq:Gaussrf} is plotted on the left in
Figure~\ref{fig:Gaussrf}. It is concave for $z\leq 1/6$ and then
convex for $z\geq 1/6$. Thus it is not possible that this LDP could
be proved, or its rate function identified, by G\"artner-Ellis style
methods that rely on convexity. The transition at $z=1/6$ occurs
when the most likely paths change from returning to $0$ within the
interval to paths that end with a non-zero position. This explains
the dramatic change in shape of the rate function at that point.

Two most likely fluid RRW paths given that $\barW_n\approx nz$ for
are shown on the right in Figure~\ref{fig:Gaussrf}. The higher path
has $z=1/3$, while the lower path has $z=1/7$. Note that, for the
lower path, all paths of this shape that start in the interval
$[0,1-\sqrt{6z}]$ are most likely paths to the deviation. That is,
there is no single most likely path, just a most likely shape that
can occur anywhere within in $[0,1]$.

{\bf M/D/1 queue-lengths}.
A second coercive example, albeit one that requires numerics for
its ultimate solution, is when $\bfmX$ is i.i.d., with $X_0$ having
Poisson distribution with rate $\alpha$ and mean $-\delta=\alpha-\mu$
($\mu\in\N$),
$P\{X_0=k\}= e^{-\alpha}\alpha^{k+\mu}/(k+\mu)!$ for
$k=-\mu,-\mu+1,\ldots$. In this setting, the RRW defined in equation
\eqref{eq:lindley} corresponds to the queue-length of an M/D/1 queue
with batch $\mu$ services, where the queue-length is observed at
customer departures. That is, at each service the minimum of
the current queue-length and $\mu$ customers are processed in a
single batch. Between services, a Poisson($\alpha$) number of
customers arrive to the queue.

The conditions of Mogul'skii's Theorem are met with
$\thetad=\thetau=+\infty$ and the local rate function is
\begin{align*}
I(x) = \alpha-(x+\mu)+(x+\mu)\log\left(\frac{x+\mu}{\alpha}\right)
\text{ if } x > -\mu
\end{align*}
and $I(x)=\infty$ if $x\le -\mu$. The sample path rate function $I_\bfmX$
is only finite at absolutely continuous functions.

Again, without loss of generality assume that $T_0=0$ and define
$T=T_1$. For a given $z>0$, for each $T\in(0,1]$ we first identify
paths $\psid$ that satisfy $\nabla I(\dot\psid(t))=\lambda^*(T-t)$ in
$[0,T]$ and $\dot\psid(T)=-\delta$, which leads to candidate optimal
paths of the following form: 
\begin{align*}
\dot\psid(t) = \alpha e^{-\lambda^*(t-T)}-\mu.
\end{align*}
Integrating, we have that
\begin{align*}
\psid(t) = 
	\frac{\alpha}{\lambda^*} e^{\lambda^*T} 
		\left(1-e^{-\lambda^*t}\right)-\mu t.
\end{align*}
If $T<1$, using the constraint $\psid(T)=0$ gives the
following equation for $\lambda^*$ in
terms of $T$
\begin{align*}
\frac{\alpha}{\lambda^*}\left(e^{\lambda^*T}-1\right)-\mu T=0.
\end{align*}
If $T=1$ and $c\eqdef \psid(1)$ we have the
following equation for $\lambda^*$ in
terms of $c$:
\begin{align*}
\frac{\alpha}{\lambda^*}\left(e^{\lambda^*}-1\right)-\mu-c=0.
\end{align*}
Both of these are transcendental equations, but can be readily
solved numerically for $\lambda^*$. Once $\lambda^*$ is known, the
constraint
\begin{align*}
z = \int_0^T \psid(t) \,dt = 
	\frac{\alpha}{\lambda^*} e^{\lambda^*T} 
		\left(T-\frac{1}{\lambda^*}\right)
		+\frac{\alpha}{(\lambda^*)^2}
		- \frac{\mu T^2}{2},
\end{align*}
gives a transcendental equation for $T$ or $c$ $(=\psid(1))$ in terms of $z$ and
identifies the solution $\psi^*$ which determines the rate
function
\begin{align*}
I_\barW(z) = \int_0^T I(\dot\psi^* (t))\,dt
	= (\alpha-\mu) T -\psi(T) 
		+\frac{\alpha}{\lambda^*}
		\left(e^{\lambda^*T}(\lambda^*T-1)+1\right).
\end{align*}

\begin{figure}[ht]
\EboxRatePath{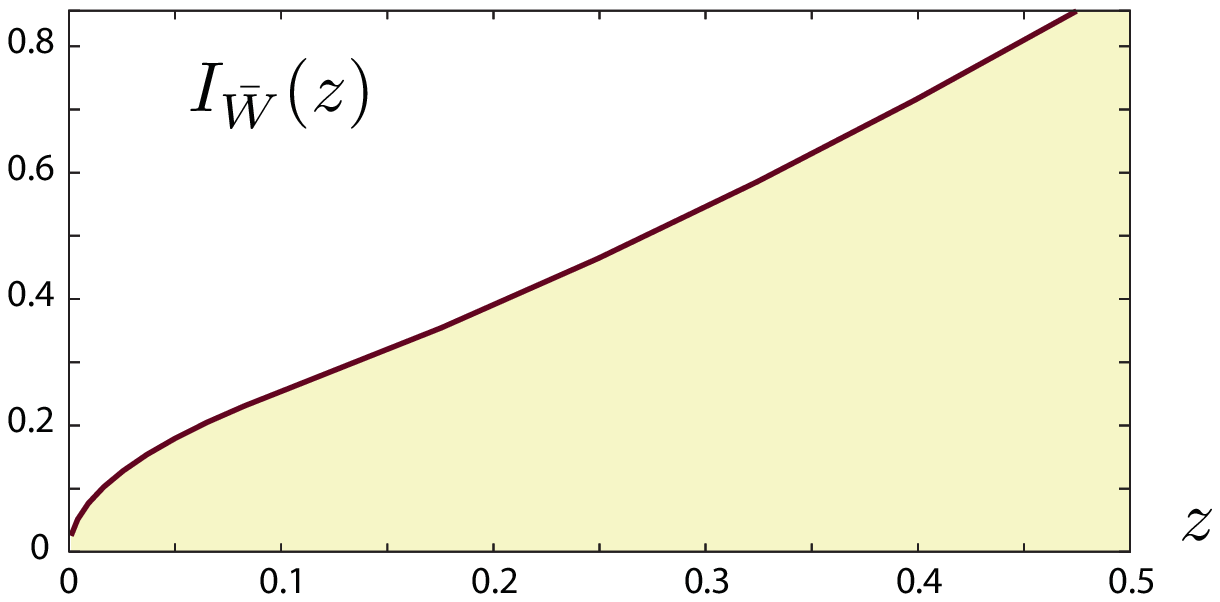}{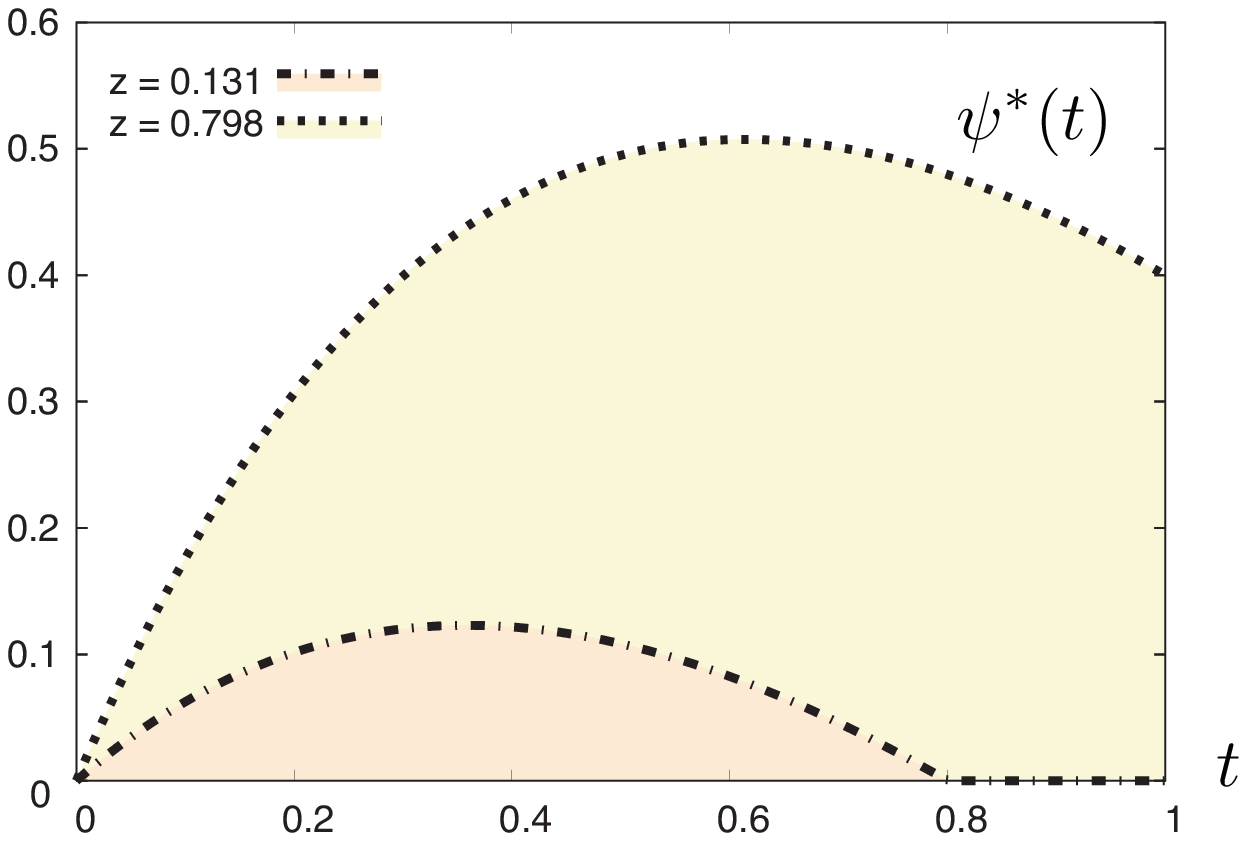}
\caption{Coercive rate function example: M/D/1 queue-lengths, Poisson
increments.
On the left hand side is $I_\barW(z)$ versus $z$. Shown on
the right hand side are most likely RRW paths,
$\phi^*(t)\approx n^{-1}W_{\lfloor nt \rfloor}$,
given that $\barW_n\approx nz$.}
\label{fig:Poissonrf}
\end{figure}

With $\alpha=0.5$ and $\mu=1.0$, the numerically calculated rate
function is plotted on the left in Figure~\ref{fig:Poissonrf}. The
transition from concave to convex again occurs when the most likely
path has both $T=1$ and $c\eqdef \psid(1)=0$. Two example most
likely paths are plotted on the right hand side of
Figure~\ref{fig:Poissonrf}, which display similar features to the
most likely paths as for the Gaussian increments RRW.

\subsection{Coercive rate function: paths with jumps}
\label{sec:examples2}

{\bf D/M/1 waiting times}.
Let $\bfmX$ be i.i.d. with $P(X_0\geq x) =\exp(-\alpha(x+\mu^{-1}))$
for $x\geq-\mu^{-1}$ giving $\E(X_0)=\alpha^{-1}-\mu^{-1}:=-\delta$.
We assume that $\mu<\alpha$ so that $\delta>0$. Then the RRW defined
in equation \eqref{eq:lindley} corresponds to waiting times at a
stable D/M/1 queue, where customers arrive at regular intervals of
length $\mu^{-1}$ and experience i.i.d. exponentially distributed
service times with mean $\alpha^{-1}$. Cram\'er's Theorem holds for
$\{x^n(1)\}$ with rate function
\begin{align*}
I(x) = 
	\begin{cases}
	\displaystyle
	\alpha\left(x+\frac 1\mu\right)
		-\log\left(\alpha\left(x+\frac 1\mu\right)\right)-1
		& \text{ if } x\in(-\mu^{-1},\infty), \\	
	\infty & \text{ otherwise},
	\end{cases}
\end{align*}
which is coercive, so that $T_0^0=T_0$ for the optimal path. However,
$\thetau=\alpha$ so that the possibility of an initial jump in the
most likely path cannot be discounted.

Again, without loss of generality, assume that $T_0=0$ and define
$T=T_1$. Using $\nabla I(\dot\psid(t))=\lambda^*(T-t)$ and the
constraint $\dot\psid(T)=-\delta=1/\alpha-1/\mu$, candidate solutions
must satisfy 
\begin{align*}
\dot\psid(t) = 
	\frac{1}{\alpha+\lambda^*(t-T)} 
	- \frac{1}{\mu}
\text{ for all } t\in(0,T)
\end{align*}
and hence, for some initial jump $\psid(0)=a\geq0$,
\begin{align*}
\psid(t) = a + 
	\frac{1}{\lambda^*}
	\log\left(\frac{\alpha+\lambda^* (t-T)}{\alpha-\lambda^*T}\right)
	- \frac{t}{\mu}
\text{ for } t\in[0,T].
\end{align*}
If $T<1$, then in addition we have that $\psid(T)=0$, which gives
the following equation:
\begin{align*}
\alpha - \lambda^* T 
	-\alpha \exp\left(\lambda^*\left(a-\frac{T}{\mu}\right)\right) = 0.
\end{align*}
If $T=1$, then we have that $c\eqdef\psid(T)\ge0$, which implies that
\begin{align*}
\alpha-\lambda^*-
	\alpha\exp\left(\lambda^*\left(a-c-\frac{1}{\mu}\right)\right)
	= 0.
\end{align*}
Treating $a$ $(=\psid(0))$ as given, these two transcendental
equations can be readily solved for $\lambda^*$.  Finally we have
the constraint that $\int_0^T\psid(t)=z$, so that
\begin{align*}
z = T\left(a-\frac{1}{\lambda^*}\right)
	+\frac{\alpha}{(\lambda^*)^2}
	\log\left(\frac{\alpha}{\alpha-\lambda^*T}\right)
	-\frac{T^2}{2\mu}
\end{align*}
and
\begin{align*}
a\alpha+\int_0^T I(\dot\psid(t))\,dt
	= \alpha\left(a+\psid(T)+\frac{T}{\mu}\right) - 2T
	+\left(\frac{\alpha-\lambda^*T}{\lambda^*}\right)
	\log\left(\frac{\alpha}{\alpha-\lambda^*T}\right).
\end{align*}
For given $a\geq0$, having solved the transcendental equations,
the most likely path and its associated rate can be calculated.    
Optimization over $a$ can then be performed numerically.

\begin{figure}[ht]
\EboxRatePath{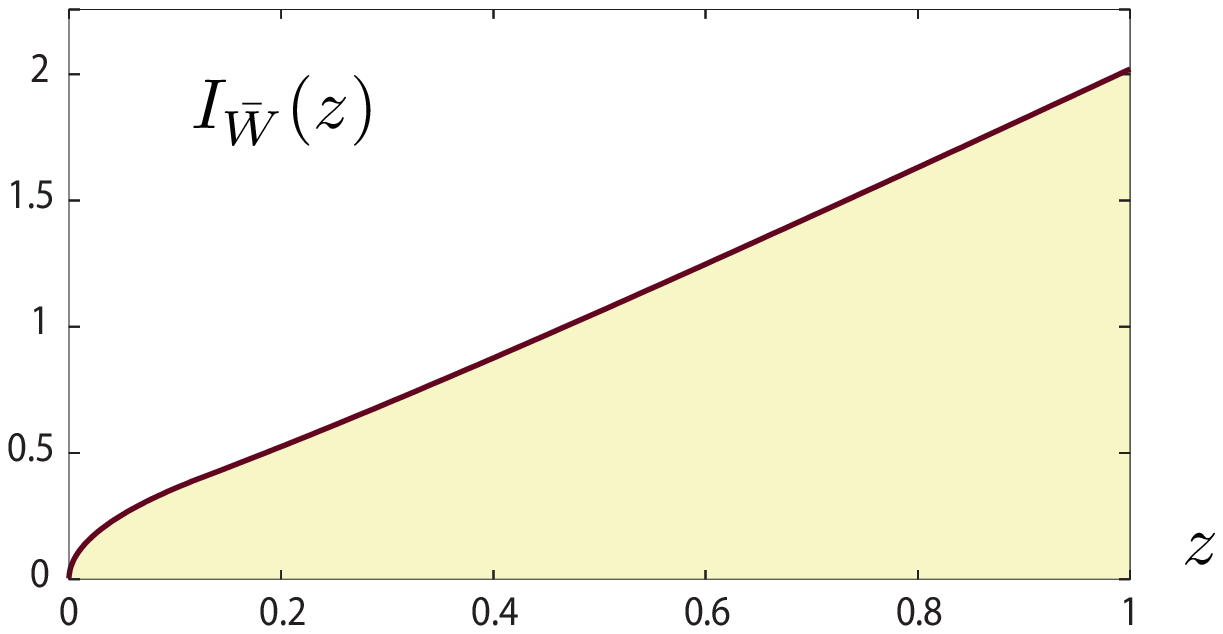}{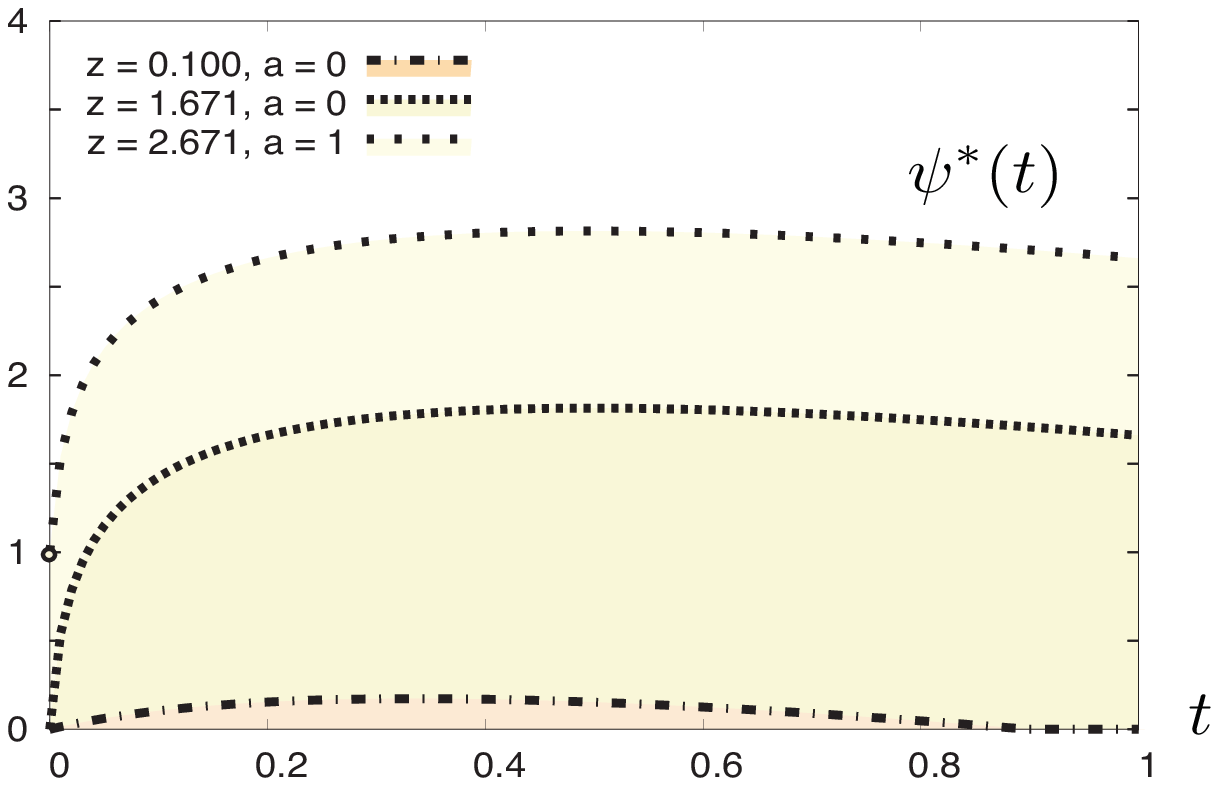}
\caption{Coercive rate function with jumps example: D/M/1 waiting times,
i.i.d. Exponential increments.
On the left hand side is $I_\barW(z)$ versus $z$. Shown on
the right hand side are most likely RRW paths,
$\phi^*(t)\approx n^{-1}W_{\lfloor nt \rfloor}$,
given that $\barW_n\approx nz$.}
\label{fig:Exponentialrf}
\end{figure}

For example, with $\alpha = 2$ and $\mu=1$, the rate function is
plotted on the left hand side of Figure~\ref{fig:Exponentialrf}.
It looks similar to the earlier examples, but is asymptotically
linear with slope $\alpha$. The reason for this is best explained
by considering the most likely path shown on the right hand side
of Figure~\ref{fig:Exponentialrf}. For small $z$, $\psid(0)=a=0$
and no jump occurs at the start of the most likely path. However,
once $z$ is sufficiently large (approximately 1.67 for these
parameters), the most likely path has a jump at $0$ followed by a
vertically shifted version of the largest-area most likely path
that doesn't have a jump, as illustrated in Figure~\ref{fig:Exponentialrf}.
The increase in the rate function for the shift of height $\psid(0)=a$
(gaining area $a$ over the interval $[0,1]$), is $a \alpha$, which
is why the rate function is ultimately linear with slope $\alpha$.

\subsection{Non-coercive rate function: rate-constrained paths}
\label{sec:examples3}

{\bf M/M/1 queue-length.} 
Let $\bfmX$ be a Bernoulli sequence taking values $-1$ and $+1$
with $\alpha=P\{X_0=+1\}<P\{X_0=-1\}=1-\alpha$. The RRW in
equation \eqref{eq:lindley} corresponds to the queue-length of an
M/M/1 queue observed at arrivals and departures.  We have $\thetau=\thetad=+\infty$.  The increments rate function $I$ is infinite outside $[-1,1]$ and is non-coercive with $\barr=1$:
\begin{align*}
I(x) = 
	\frac{1+x}{2}\log\left(\frac{1+x}{2\alpha}\right)
	+\frac{1-x}{2}\log\left(\frac{1-x}{2(1-\alpha)}\right)
\text{ and } \barr=1.
\end{align*}

Note that $I_\barW(z)=+\infty$ if $z>1/2$; if $z=1/2$, then $T_0^0=T=1$
and $\psi^*(t)=t$ so that $I_\barW(z)= -\log(\alpha)$.
Without loss of generality, let $T_0=0$ and define $T=T_1$.
Assume that $T_0^0<T$. The equation $\nabla
I(\dot\psid(t))=\lambda^*(T-t)$ with the boundary condition 
$\dot\psid(T)=-\delta =2\alpha-1$ gives
\begin{align*}
\dot\psid(t) = 
	\begin{cases}
	1 & \text{ if } t\in[0,T_0^0] \\
	\displaystyle
	\frac{\alpha\exp(2\lambda^*(T-t))-(1-\alpha)}
		{\alpha\exp(2\lambda^*(T-t))+1-\alpha}
		& \text{ if } t\in(T_0^0,T]. 
	\end{cases}
\end{align*}
By integrating, we are looking at proposed solutions
\begin{align*}
\psid(t) = 
	\begin{cases}
	t	& \text{ if } t\in[0,T_0^0] \\
	\displaystyle
		2T_0^0-t
		+\frac{1}{\lambda^*}
		\log\left(
		\frac
		{\alpha e^{2\lambda^*(T-T_0^0)}+1-\alpha}	
		{\alpha e^{2\lambda^*(T-t)}+1-\alpha}	
			\right)
		& \text{ if } t\in(T_0^0,T]. 
	\end{cases}
\end{align*}
If $T<1$, then $\psid(T)=0$ gives
\begin{align*}
\alpha e^{2\lambda^*(T-T_0^0)} - e^{\lambda^*(T-2T_0^0)} +1-\alpha=0
\end{align*}
and, in particular, if $T_0^0=0$, then
$\lambda^*=\log((1-\alpha)/\alpha)/T$.  While if $T=1$, then
$c\eqdef\psid(T)\ge 0$ ($c<1$) gives the following equation for
$\lambda^*$
\begin{align}
\label{eq:MM1lambda}
\alpha e^{2\lambda^*(1-T_0^0)} - e^{\lambda^*(c+1-2T_0^0)} +1-\alpha=0.
\end{align}
Note that this equation only has a positive solution for 
$c\in(\max\{0,2\alpha-1+2(1-\alpha)T_0^0\},1)$. The lower bound
embodies the fact that the optimal path has $\dot\psi^*(t)\geq
2\alpha-1$ and therefore $c=\psi^*(1)=T_0^0+\int_{T_0^0}^1\dot\psi^*(t)\,dt
\geq 2\alpha-1+2(1-\alpha)T_0^0$.
Once $\lambda^*$ is identified, one can numerically evaluate the
integral
\begin{align*}
z = \int_0^T \psid(t)\, dt.
\end{align*}

\begin{figure}[h]
\EboxRatePath{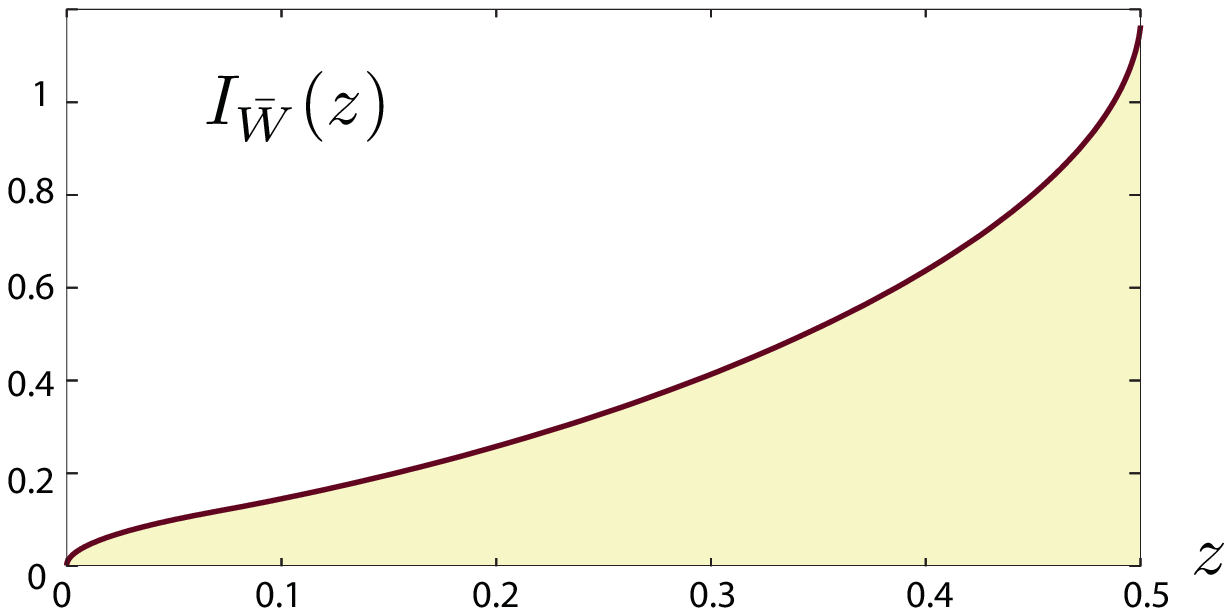}{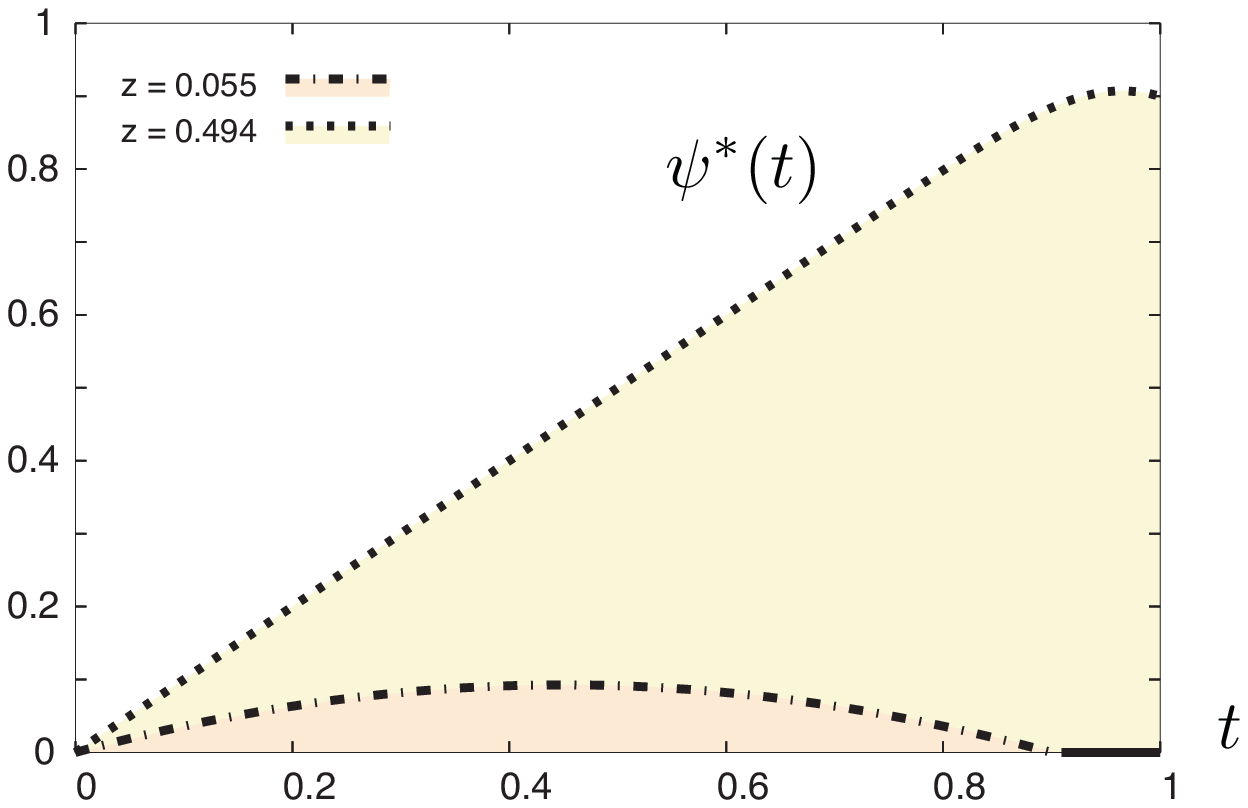}
\caption{Non-coercive rate function example: M/M/1 queue-lengths, Bernoulli
increments.
On the left hand side is $I_\barW(z)$ versus $z$. Shown on
the right hand side are most likely RRW paths,
$\phi^*(t)\approx n^{-1}W_{\lfloor nt \rfloor}$,
given that $\barW_n\approx nz$.}
\label{fig:Bernoullirf}
\end{figure}

For example, with $\alpha=1/3$, Figure~\ref{fig:Bernoullirf} plots
the numerically evaluated rate function $I_\barW(z)$ versus $z$.
The initial shape of the rate function is similar to $\sqrt{x}$,
but, as can be seen clearly in the graph, once $z$ is sufficiently
large that the optimal path has $c>0$, the rate function increases
dramatically.

For $z<0.0682$, the optimal value of $T_1$ is less than $1$ and the
numerically-identified most likely path occurs with $T_0^0=0$.  For
$z>0.0682$, $T=1$ and the optimal path also has $T_0^0=0$ apart from,
possibly, as $z\to 1/2$. The reason for this caveat is that the
cost of a path with $T_0^0=0$ becomes numerically indistinguishable
from those with $T_0^0>0$ if $z\approx1/2$. To see this, note that
as $z\uparrow 1/2$, $c\uparrow1$, so that if $T_0^0=0$, then
$\lambda^*$, the solution of equation \eqref{eq:MM1lambda}, tends
to $1$ and the most likely path to large simulated mean has slope
close to $1$ for a substantial range of $t$. This is illustrated
in the higher paths on the right in Figure~\ref{fig:Bernoullirf}
corresponding to the most likely path for the deviation $z\approx0.45$.
For this path $\lambda^*\approx0.941$. Note that the slope is nearly
$1$ until near $t=0.85$, even though technically $T_0^0=0$.

\subsection{Simulations}

One of the strong deductions of these sample-path arguments is the
prediction of the most likely path that gives rise to a large
simulated mean. To illustrate the merits of these predictions we
conducted simulations of the RRW in two settings: Gaussian i.i.d.\
increments, and the M/M/1 queue example introduced in
Section~\ref{sec:examples3}.

\begin{figure}[h]
\begin{center}
\parbox[t]{.425\hsize}{\epsfxsize=\hsize \epsfbox{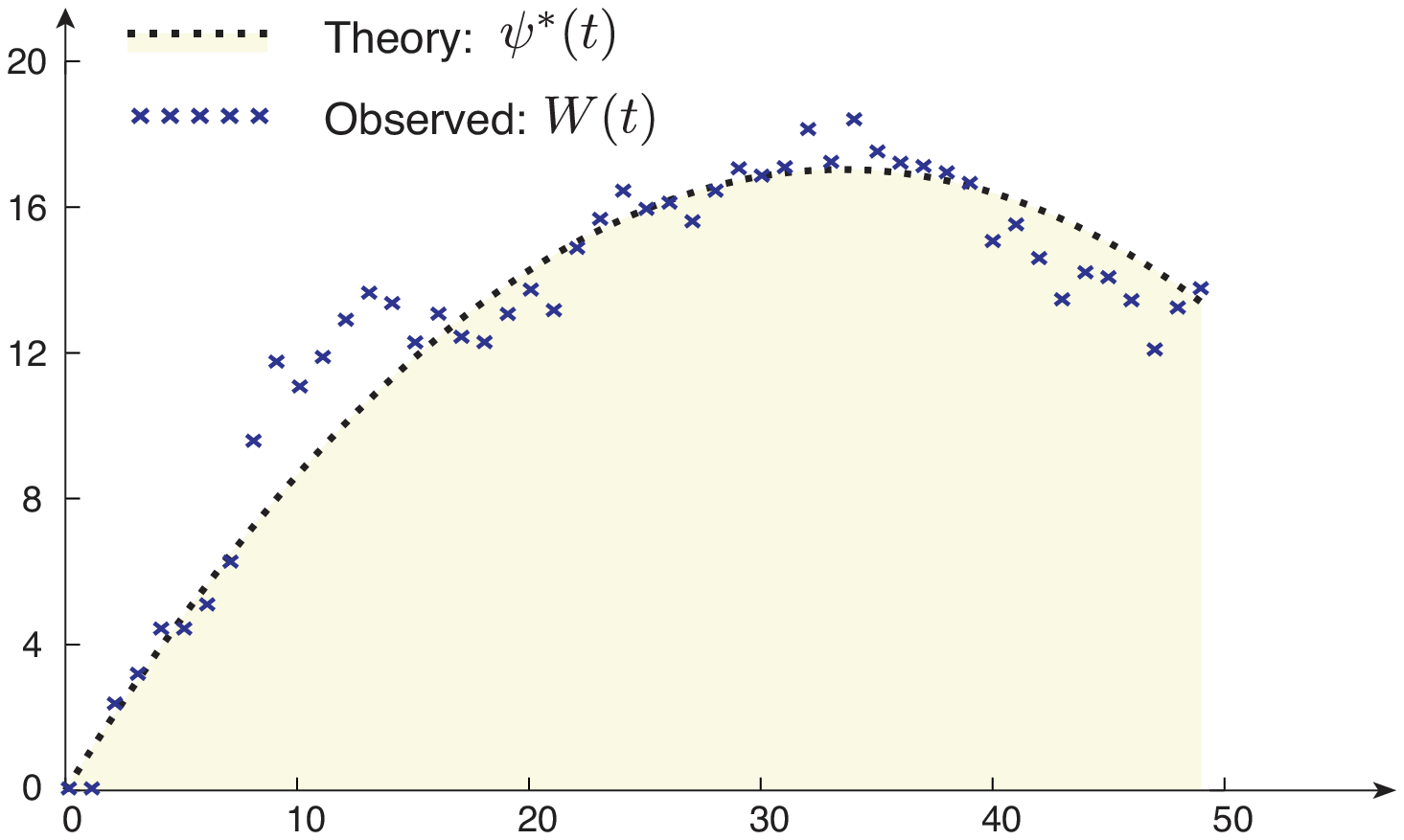}}
\hfil
\parbox[t]{.425\hsize}{\epsfxsize=\hsize \epsfbox{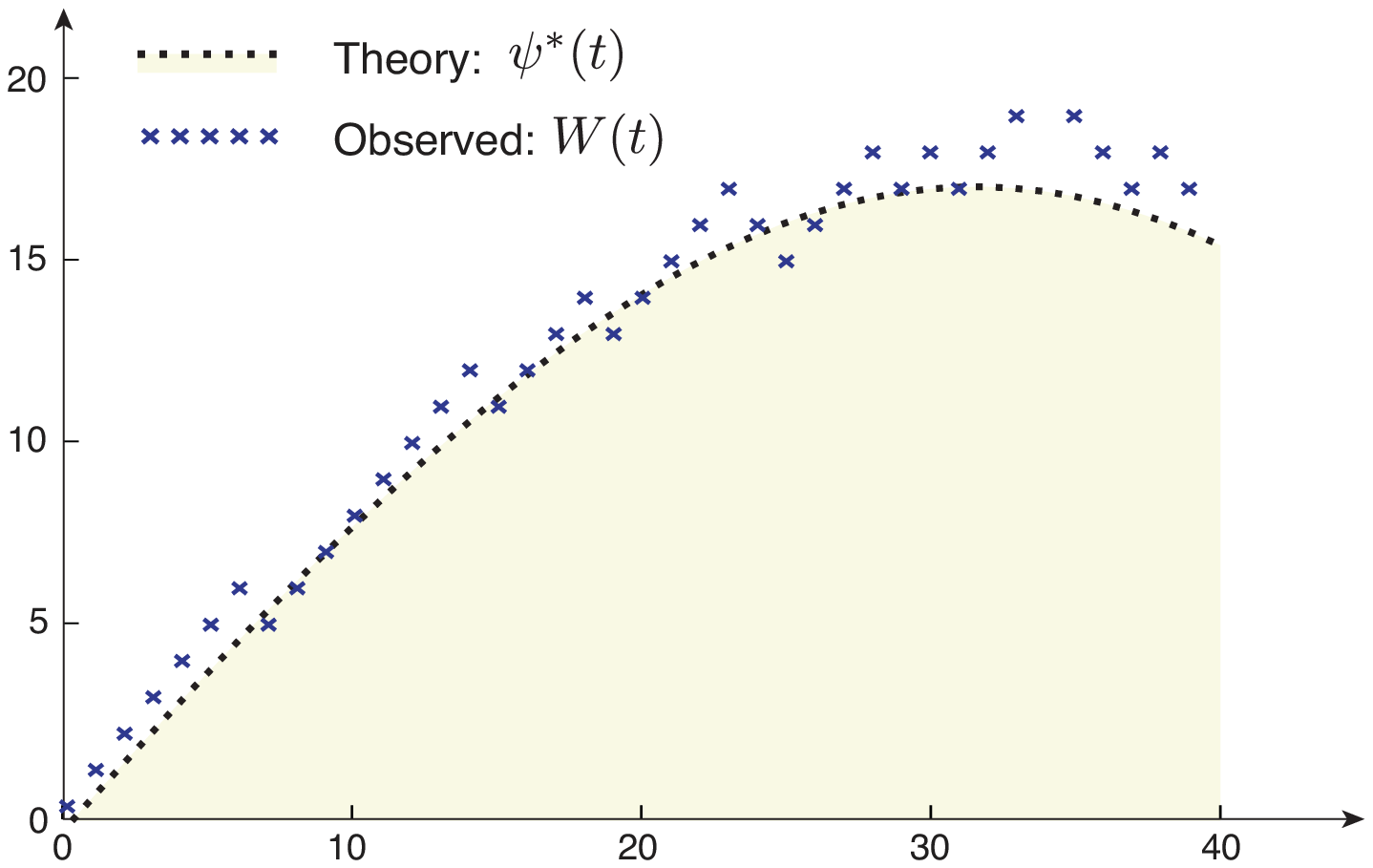}}
\end{center}
\caption{RRW with i.i.d.\ increments. The figure shown on the left
hand side shows experiments obtained with Gaussian increments. On
the right hand side are results obtained for the M/M/1 queue model
in which the increments take on values $\pm 1$. In each case,
the observed path has the largest simulated mean out of $2\times
10^8$ sampled paths. Also, shown in each figure is the corresponding
theoretical prediction of the most likely path, given the observed
simulated mean.}
\label{fig:GaussianMM1simulated}
\end{figure}

In each case, the RRW was simulated for a fixed number of steps
$n\ge 1$, and the simulation was repeated $2\times10^8$ times. Of
these simulated RRWs, the one with the largest simulated mean was
recorded and compared with the theory laid out in
Section~\ref{sec:examples}. This theory predicts the approximation
$W_{[t]}\approx n\,\phi^*(t/n)$ for $t\in[0,n]$. The results from
two experiments are illustrated in Figure~\ref{fig:GaussianMM1simulated}.

In the first experiment illustrated on the left hand side, the
increments of the random walk were taken to be i.i.d. Gaussian with
$\delta=0.5$, $\sigma^2=1$ and the time-horizon $n=50$. The second
experiment used the M/M/1 queue example found in
Section~\ref{sec:examples3} with $\alpha=0.3$ and $n=40$. In each
experiment, the observed sample path is plotted along with the
theoretical prediction corresponding to the observed simulated mean.
The theory's quantitative power in predicting the shape of the most
likely path is apparent in these figures.

\section{Discussion}

As a final remark, we have mentioned that our fundamental hypothesis,
Assumption \ref{ass:LDPX}, encompasses the light tailed setting in
the absence of long range dependence. However, by changing the speed
of the LDP, this assumption also holds for certain long range
dependent processes.  For example, in continuous time it is known,
e.g. \cite{Deuschel89}\cite{Norros99}\cite{Majewski03}, that
fractional Brownian Motion (fBM) with Hurst parameter $H$ satisfies
the LDP at speed $n^{2(1-H)}$ in $D[0,1]$. As the nature of the
speed does not enter the proof of Theorem \ref{thm:waiting}, the
first conclusion of that theorem holds for these processes. However,
even for the canonical example of fBM the rate function is not of
the integral form in equation \eqref{eq:IDZ} and thus, in the long
range dependent setting, it is hard to deduce if any general
properties exist for the most likely paths to a large simulated
mean.

\paragraph{Acknowledgment}
{\footnotesize
Financial support from the National Science Foundation (NSF CCF
07-29031), AFOSR (FA9550-09-1-0190) and Science Foundation Ireland
(07/IN.1/1901) is gratefully acknowledged.  Any opinions, findings,
and conclusions or recommendations expressed in this material are
those of the authors and do not necessarily reflect the views of
NSF, AFOSR or SFI.}

\def\cprime{$'$} \def\cprime{$'$} \def\cprime{$'$}
\providecommand{\bysame}{\leavevmode\hbox to3em{\hrulefill}\thinspace}
\providecommand{\MR}{\relax\ifhmode\unskip\space\fi MR }
\providecommand{\MRhref}[2]{%
  \href{http://www.ams.org/mathscinet-getitem?mr=#1}{#2}
}
\providecommand{\href}[2]{#2}

\end{document}